\documentstyle [12pt,amssymb,psfig]{article}
\title{Some free-by-cyclic groups}
\author{Ian J. Leary \thanks{Partially supported by EPSRC grant 
no.\ GR/L69398} 
\and Graham A. Niblo \thanks{Partially supported by EPSRC grant 
no.\ GR/K25618}
\and Daniel T. Wise \thanks{Supported by NSF grant no.\ DMS-9627506} }
\def\figone{\vskip 7mm\centerline{\psfig{figure=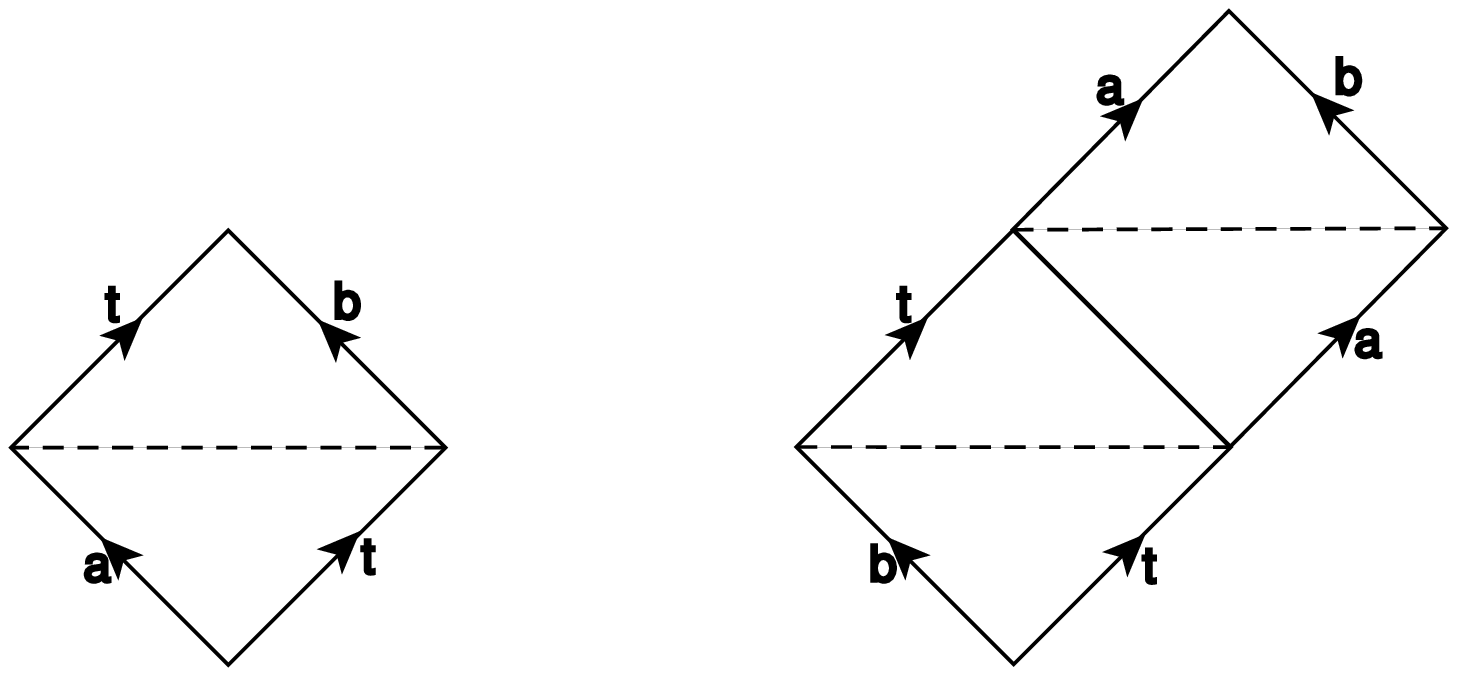,height=2in}}
\centerline{Figure 1.}\vskip 7mm}
\def\figtwo{\vskip
7mm\centerline{\vbox{\hbox{\psfig{figure=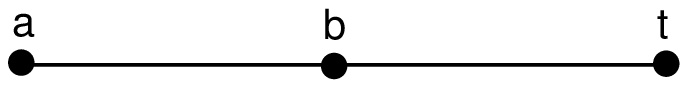,height=.25in}}
\hbox{\quad\ \  The descending link}}
\qquad\vbox{\hbox{\psfig{figure=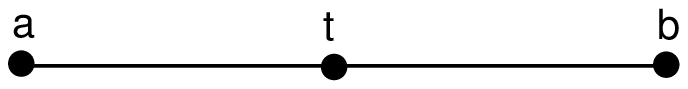,height=.25in}}
\hbox{\quad\ \  The ascending link}}}
\centerline{Figure 2.}\vskip 7mm}
\def\figthree{\vskip 7mm\centerline{\psfig{figure=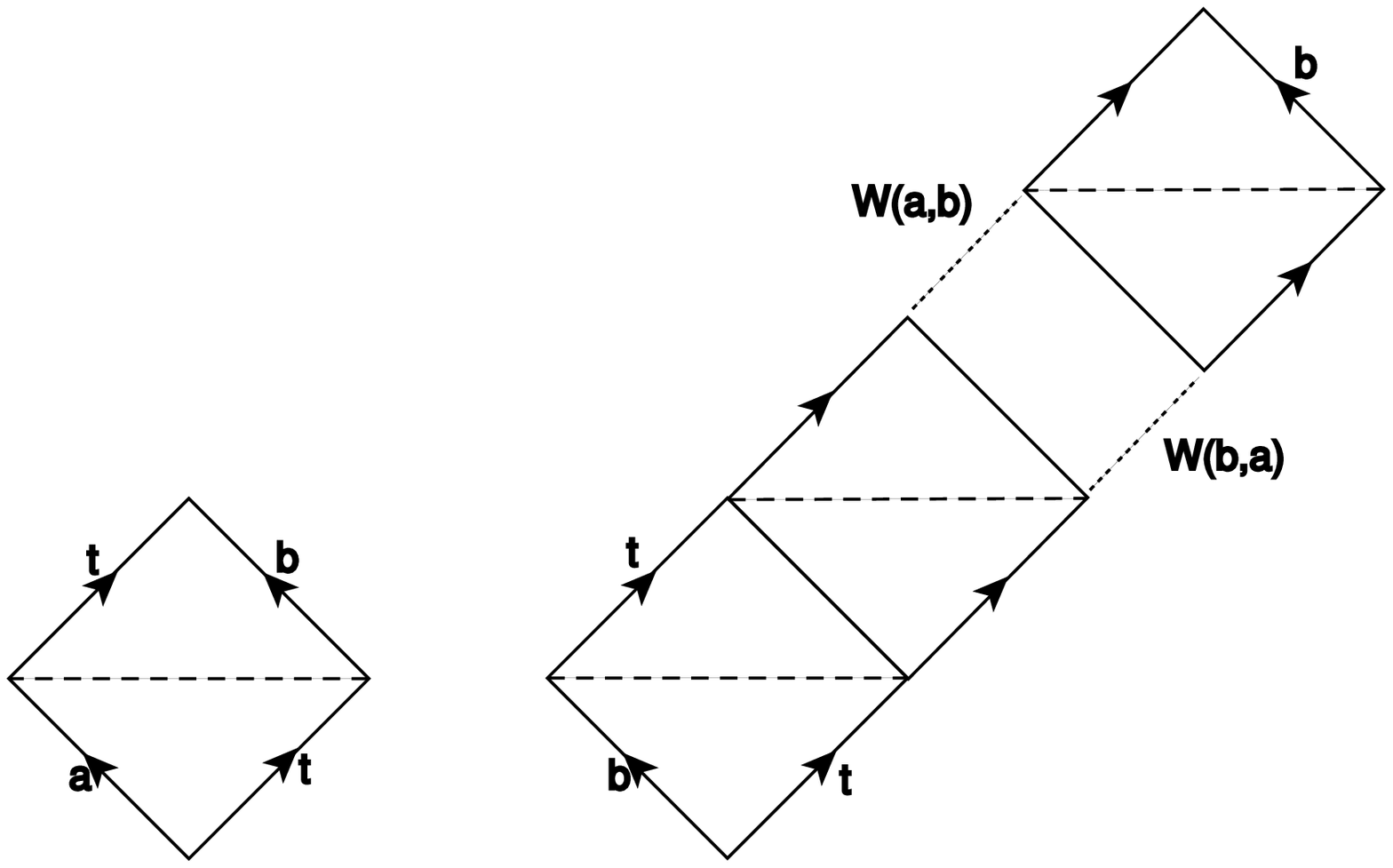,height=3in}}
\nobreak\centerline{Figure 3.}\vskip 7mm}

\begin{document}
\maketitle
\noindent
A group is said to be locally free if every finitely generated
subgroup of it is free.  One example is the additive group of the
rationals.  We exhibit a finitely generated group $G$ that is
free-by-cyclic and contains a non-free, locally free subgroup. 
The smallest such example that we have found is of the form
$G\cong F_n\rtimes\Bbb Z$ for $n=3$.  We also construct
word-hyperbolic 
examples for larger values of $n$, and show that the groups are
not subgroup separable.  We used Bestvina and Brady's
\lq Morse theory for cube complexes' in the construction of these
groups.  The authors thank Jim Anderson, who posed a question
concerning 3-manifolds that led to these examples, and the 
referee, whose comments were very helpful.   
This work was started at a conference at Southampton, immediately 
before Groups St.\ Andrews, which was funded by EPSRC visitor grants 
GR/L06928 and GR/L31135, and by a grant from the LMS.

Throughout this note, $F_n$ denotes a free group of rank $n$, 
$\bar x$ denotes $x^{-1}$, and $x^y = \bar y x y$.

\proclaim Proposition 1. The group $G$ given by the presentation
$$G=\left\langle a,b,t\, \colon a^t= b, \,\,
b^t = ab\bar a\right\rangle$$
contains a non-free, locally free subgroup, and is isomorphic to
a split extension $F_3\rtimes \Bbb Z$. 
\par 
\noindent
{\bf Proof } The given presentation expresses $G$ as an ascending
HNN extension, with base group freely generated by $a$ and $b$ and
with stable letter $t$.  Define $\phi \colon G\rightarrow \Bbb Z$ by
$\phi(t)=1,\, \phi(a)=\phi(b)= 0$, and let $K$ be the kernel of
$\phi$.  Then $K$ is a strictly ascending union of 2-generator free
groups
$$\langle a,b\rangle \subseteq  \langle a,b \rangle ^{\bar t}
\subseteq  \langle a,b \rangle ^{\bar t^2} \subseteq \cdots
\subseteq K.$$
Any such group is locally free and not free, since it is not
finitely
generated and the rank of its abelianization is at most two.  (In
fact, the abelianization of $K$ is infinite cyclic.) 

It remains to show that $G$ is free-by-cyclic.  Define $\psi \colon
G\rightarrow \Bbb Z$ by $\psi(t)=\psi(a)=\psi(b)=1$.  It will be shown
that the kernel of $\psi$ is free of rank 3.  A presentation 2-
complex
$Y$ for $G$ may be constructed by attaching two 2-cells to a rose
with
edges $a$, $b$, and $t$ according to the maps given in figure~1.
Since $Y$ is obtained from the presentation of $G$ as an
HNN extension with free base group, it follows that $Y$ is
an Eilenberg-Mac~Lane space for $G$ (see proposition~3.6 of
\cite{scw}).  Represent the three 1-cells of $Y$ as unit intervals, 
and represent the two 2-cells of $Y$ as a unit
square and a $2\times 1$ rectangle, as indicated in figure~1. 
This makes $Y$ into an affine cell complex in the sense of Bestvina 
and Brady (\cite{bb}, Def.\ 2.1). 

\figone

Now take $S^1={\Bbb R}/{\Bbb Z}$, viewed as a cell complex with one
vertex and one edge of length 1, as an Eilenberg-Mac~Lane space
for the integers.  A cellular map $g\colon Y \rightarrow S^1$
may be defined that induces the homomorphism
$\psi\colon G\rightarrow \Bbb Z$
on fundamental groups and is affine on each cell.  In figure~1
this map is represented by \lq height modulo one', where the length 
of each edge is chosen so that its height is one.  The inverse 
image of the vertex $v$ of $S^1$ is a rose consisting of one vertex 
and three 1-cells (the dotted
lines on figure~1).  Now let $X$ be the cover of $Y$ corresponding
to the subgroup $H=\ker(\psi)$.  The map $g$ lifts to a map $f\colon
X\rightarrow \Bbb R$.  $X$ is an affine cell complex, and $f$ is a
Morse function in the sense of \cite{bb}, Def.\ 2.2.  By construction,
$X$ is an Eilenberg-Mac~Lane space for $H$, and for any integer $t$,
$X_t=f^{-1}(t)$ consists of a disjoint union of copies of a 
3-petalled rose.  (Clearly, $X_t$ is a disjoint union of connected 
covers of $g^{-1}(v)$, but since any loop in $g^{-1}(v)$ 
represents an element of $G=\pi_1(Y)$ in the kernel of $\psi$, every
lift in $X_t$ of such a loop is itself a loop, and hence $X_t$ is a 
disjoint union of 1-fold covers.)  

Bestvina and Brady's Morse theory allows one to compare $X$ and
$X_t$: by lemma~2.5 of \cite{bb}, a space homotopy equivalent to $X$
may be obtained from $X_t$ by coning off a subspace homeomorphic
to a copy of the descending link
(resp.\ ascending link) at $v$ for each vertex $v$ of $X$ such that
$f(v)>t$ (resp.\ $f(v)<t$).  (Ascending and descending links are
defined in section~2 of \cite{bb}.)  All vertices of $X$ have isomorphic
links, since $Y$ has only one vertex, and the ascending and
descending links at each vertex are as shown in figure~2. 
\figtwo
\noindent
Both the ascending and 
descending link are contractible.  Since coning off a contractible
subspace does not change the homotopy type of a space, it follows
that $X$ is homotopy equivalent to $X_t$.  But it is already known that
$X$ is an Eilenberg-Mac~Lane space for $H$, and that $X_t$ is a
disjoint union of 3-petalled roses.  It follows that $X_t$ is
connected, and that $H$ is free of rank three. 
\quad\hbox{\vrule width 3pt height 6pt depth 1pt}

With the benefit of hindsight, a shorter proof that $G$ as above is
free-by-cyclic may be given---see Proposition~2 below.  Such a proof
gives no indication as to how $G$ was discovered however.  Moreover,
the techniques of Proposition~1 generalize easily to more complicated
presentations such as those given in Proposition~3.

\proclaim Proposition 2.  Let $H'$ be freely generated by $x$, $y$
and $z$, and define an automorphism $\theta$ of $H'$ by
$$\theta(x)=y,\quad \theta(y)=z, \quad \theta(z)= y^2\bar x.$$
The group $G$ of Proposition~1 is isomorphic to $H'\rtimes\langle
t\rangle$, where the conjugation action of $t$ on $H'$ is given by
$\theta$. 
\par 
\noindent
{\bf Proof } First, check that the endomorphism $\theta$ is an
automorphism of $H'$ by exhibiting an inverse:
$$\theta^{-1}(z)=y,\quad \theta^{-1}(y)=x, \quad
\theta^{-1}(x)= \bar z x^2.$$
Now, eliminate $b$ from the given presentation for $G$ to
obtain
$$G=\langle a,t\,\colon \bar t \bar t a tt a \bar t \bar a t \bar
a\rangle.$$
Substitute $a=xt$, and eliminate $a$, obtaining
$$G=\langle x,t\,\colon \bar t \bar t (xt) tt (xt) \bar t
(\bar t\bar x)  t (\bar t \bar x) \rangle =
\langle x,t\,\colon \bar t^2 xt^3 x\bar t \bar x^2\rangle.$$
Add new generators $y=x^t$ and $z= x^{t^2}$, obtaining
\begin{eqnarray}
G&=& \langle x,y,z,t\,\colon x^t=y,\, y^t=z,\,
ztx\bar y^2 \bar t\rangle\nonumber\\
 &=& \langle x,y,z,t\, \colon
x^t= y,\, y^t = z, \, z^t = y^2\bar x \rangle.\nonumber
\end{eqnarray}
Thus $G$ is seen to be isomorphic to $H'\rtimes \langle t \rangle$ as
claimed.  \qquad
\hbox{\vrule width 3pt height 6pt depth 1pt}

Next, we show how to construct a word-hyperbolic group having
similar properties to the group $G$. 

\proclaim Proposition 3.  For $s\geq 3$ define a word $W(x,y)$ by
$$W(x,y)= xy^4xy^5x\cdots xy^{4+s}x.$$ 
The group $G_s$ with presentation
$$G_s = \langle a,b,t\,\colon a^t=b,\, b^t = W(b,a)b(W(a,b))^{-1}
\rangle$$
is free-by-cyclic and contains a non-free, locally free subgroup. 
For
$s$ sufficiently large, $G_s$ is word-hyperbolic. 
\par
\noindent
{\bf Proof } As in Proposition~1, $G_s$ is a strictly ascending
HNN-extension with base group freely generated by $a$ and $b$, so
contains a non-free, locally free subgroup.  As in Proposition~1,
an
Eilenberg-Mac~Lane space for $G_s$ with an affine cell structure
can
be made by attaching a unit square and an $m\times 1$ rectangle to
a rose with three petals of length 1.  (Here $m$ is one more than
the length of the word $W$, as shown in figure~3.)  The argument
given
in Proposition~1 shows that $G_s$ is expressible as $F_n\rtimes \Bbb
Z$,
where $n$ is the total area of the two 2-cells in figure~3, i.e.,
$n = m+1 = s+ 4 + (8+s)(s+1)/2$. 

\figthree

It remains to show that $G_s$ is word-hyperbolic for $s$
sufficiently
large.  For this, it suffices to show that some presentation
for $G_s$ satisfies the $C'(1/7)$
small cancellation condition (see \cite{str}).  Eliminate $b$ from the
presentation
for $G_s$.  This leaves a 1-relator group, with relator
$$bt a b^4 a b^5 a \cdots a b^{4+s}a \bar b^2 \bar a^{4+s} \bar
b \cdots \bar b \bar a^5 \bar b \bar a^4 \bar b \bar t\qquad$$
$$\quad = \bar t a t^2 a (\bar t a^4 t a )(\bar t a^5 t a) 
\cdots (\bar t a^{4+s} ta) \bar t \bar a^2 (t \bar a^{4+s} \bar t
\bar a)( t \bar a^{3+s} \bar t \bar a)\cdots  (t \bar a^5 \bar t
\bar a) (t \bar a^4\bar t \bar a).$$
The total length of this relator as a cyclic word 
in $a$ and $t$ may be seen
to be $8 + (14+s)(s+1)$.  (The bracketing of the word is intended to
facilitate this check.)  For any $4\leq r \leq 4+s$, the four words 
$(\bar t a^r t)^{\pm 1}$ and $(t a^r \bar t)^{\pm 1}$ occur exactly
once each as a subword of the relator or its inverse, and any subword
of the relator or its inverse of length at least $2s+15$ contains a
subword of this form.  (The worst case is the subword 
$a^{4+s}ta\bar t \bar a^2 t \bar a^{4+s}$, 
of length $2s+14$.)  Hence any subword of the
relator or its inverse of length $2s+15$ occurs in a unique place. 
It follows that $G_s$ is word-hyperbolic whenever
$2s + 15 \leq 1/7\left( 8 + (14+s)(s+1)\right)$.  This inequality is
satisfied for all sufficiently large~$s$.  (In fact, $s\geq 9$
suffices.)  \qquad \hbox{\vrule width 3pt height 6pt depth 1pt}

P. Scott asked if free-by-cyclic groups are necessarily subgroup
separable.  An example due to Burns, Karass and Solitar showed that
this is not the case (see \cite{bks}).  The groups constructed above give
another, simpler, argument to show this. 

\proclaim Proposition 4.  The groups $G$ and $G_s$, constructed in
Propositions 1~and~3, are not subgroup separable. 
\par
\noindent
{\bf Proof } In each case, let $L_1$ be the subgroup generated by $a$
and $b$, and let $L_2 = tL_1\bar t$.  Then $L_1$ and $L_2$ are free of
rank two and are conjugate in $G$ (resp.\ in $G_s$).  Moreover, $L_1$
is a proper subgroup of $L_2$.  It follows that $L_1$ cannot be
closed, since it cannot be separated from any element of
$L_2\setminus
L_1$:  in any finite quotient, the images of $L_1$ and $L_2$ have the
same order, since they are conjugate, and so must be equal since the
image of $L_1$ is a subgroup of the image of $L_2$. 
\qquad \hbox{\vrule width 3pt height 6pt depth 1pt}

\vfill\eject

\vskip .5in
\obeylines
Authors' addresses: 
Faculty of Mathematical Studies,
University of Southampton,
Southampton,
SO17 1BJ. 
\vskip .23456789012345678901in
Department of Mathematics,
Cornell University,
Ithaca, NY
14853

\end{document}